\documentclass[11pt]{article}
\usepackage{bbm}
\usepackage{mathrsfs}
\usepackage{amscd}
\usepackage{amsmath,amsfonts,amssymb,amscd}
\usepackage{indentfirst,graphics,epsfig,psfrag}
\input{epsf}
\usepackage{ifpdf}
\usepackage{enumerate}
\usepackage{appendix}
\usepackage{enumerate}

\usepackage{lineno}

\makeatletter \@addtoreset{figure}{section} \makeatother
\makeatletter
\long\def\@makecaption#1#2{%
   \vskip 10\p@
   \setbox\@tempboxa\hbox{{#1}\ \ #2}%
   \ifdim \wd\@tempboxa >\hsize

       {#1}\ \ #2\par
   \else
       \hbox to\hsize{\hfil\box\@tempboxa\hfil}%
   \fi}
\makeatother

\newtheorem{thm}{Theorem}
\newtheorem{cor}{Corollary}
\newtheorem{lem}{Lemma}

\newtheorem{obs}{Observation}
\newtheorem{pro}{Proposition}

\newcommand{\qed}{{\hfill\rule{3pt}{7pt}}}

\def\qed{\hfill \rule{4pt}{7pt}}

\begin{document}
\title{\textbf{Monochromatic connectivity and graph products} \footnote{Supported by the National Science Foundation of China
(Nos. 11161037, 11101232, 11461054) and the Science Found of Qinghai Province
(No. 2014-ZJ-907).}}
\author{
\small Yaping Mao, \ Zhao Wang \footnote{Corresponding author},  \ Fengnan Yanling, \ Chengfu Ye\\[0.2cm]
\small Department of Mathematics, Qinghai Normal\\[-1mm]
\small University, Xining, Qinghai 810008, China\\[0.2cm]
\small E-mails: maoyaping@ymail.com; wangzhao380@yahoo.com;\\
\small fengnanyanlin@yahoo.com; yechf@qhnu.edu.cn}
\date{}
\maketitle
\begin{abstract}
The concept of monochromatic connectivity was introduced by Caro and Yuster.
A path in an edge-colored graph is called a \emph{monochromatic path} if all the edges on
the path are colored the same. An edge-coloring of $G$ is a \emph{monochromatic connection
coloring} ($MC$-coloring, for short) if there is a monochromatic path joining any two
vertices in $G$. The \emph{monochromatic connection number}, denoted by $mc(G)$, is defined
to be the maximum number of colors used in an $MC$-coloring of a graph $G$. In this paper, we study the monochromatic connection number on the lexicographical, strong, Cartesian
and direct product and present several upper and lower bounds for these
products of graphs.\\[2mm]
\end{abstract}

\section{Introduction}

All graphs considered in this paper are simple, finite and undirected. We follow the
terminology and notation of Bondy and Murty \cite{Bondy}. For a graph $G$, we use $V(G)$, $E(G)$,
$n(G)$, $m(G)$, $\delta(G)$, $\kappa(G)$, $\kappa'(G)$, $\delta(G)$ and $diam(G)$ to denote the vertex set, edge set, number of vertices, number of edges, connectivity, edge-connectivity, minimum degree and diameter of $G$, respectively. The rainbow connections of a graph which are applied to measure the safety of a network are introduced by Chartrand, Johns, McKeon and
Zhang \cite{Chartrand}. Readers can see \cite{Chartrand, ChartrandJMZ,
ChartrandOZ} for details. Consider an edge-coloring (not necessarily
proper) of a graph $G=(V,E)$. We say that a path of $G$ is
\emph{rainbow}, if no two edges on the path have the same color. An
edge-colored graph $G$ is \emph{rainbow connected} if every two
vertices are connected by a rainbow path. The minimum number of
colors required to rainbow color a graph $G$ is called \emph{the
rainbow connection number}, denoted by $rc(G)$. For more
results on the rainbow connection, we
refer to the survey paper \cite{LiSS} of Li, Shi and Sun and a new
book \cite{LiS} of Li and Sun.

Let $G$ be a nontrivial connected graph with an \emph{edge-coloring} $f: E(G)\rightarrow \{1,2,\ldots ,\ell\}$,
$\ell\in N$, where adjacent edges may be colored the same. A path of $G$ is a \emph{monochromatic path} if all the edges on the path are colored the same. An edge-coloring of $G$ is a \emph{monochromatic connection coloring} ($MC$-coloring, for short) if there is a monochromatic path joining any two vertices in $G$. How colorful can an $MC$-coloring be ? One can see that this question
is the natural opposite of the well-studied problem on rainbow connection number of graphs. Let $mc(G)$ denote the maximum number of colors used in an $MC$-coloring of a graph $G$, which called the \emph{monochromatic connection number} of $G$. Note that an $MC$-coloring does not exist if G is not connected, and in this case we simply let $mc(G)=0$.

These concepts were introduced by Caro and Yuster in \cite{CaroY}. For more results on monochromatic connection number, we refer to \cite{CaiLW, CaiLW2, CaroY, GuLQ}. The following observation is immediate.

\begin{obs}{\upshape \cite{CaroY}}\label{obs1}
Let $G$ be a connected graph with $n(G)$ vertices and $m(G)$ edges. Then
$$
mc(G)\geq m(G)-n(G)+2.
$$
\end{obs}

Simply color the edges of a spanning tree with one
color, and each of the remaining edges may be assigned a distinct fresh color. Caro and
Yuster gave some sufficient conditions for graphs attaining this lower bound.

\begin{thm}{\upshape \cite{CaroY}}\label{th1}
Let $G$ be a connected graph with $n>3$. If $G$ satisfies any of the
following properties, then $mc(G)=m-n+2$.

$(a)$ $G$ (the complement of $G$) is $4$-connected;

$(b)$ $G$ is triangle-free;

$(c)$ $\Delta(G)<n-\frac{2m-3(n-1)}{n-3}$; In particular, this holds if $\Delta(G)\leq (n+1)/2$, and this also
holds if $\Delta(G)\leq n-2m/n$.

$(d)$ $Diam(G)\geq 3$;

$(e)$ $G$ has a cut vertex.
\end{thm}

Product networks were proposed based upon the idea of using the
cross product as a tool for ``combining'' two known graphs with
established properties to obtain a new one that inherits properties
from both \cite{DayA}. Recently, there has been an increasing
interest in a class of interconnection networks called Cartesian
product networks; see \cite{Bao, DayA}.
The other standard products (Direct, strong, and lexicographic)
draw a constant attention of graph research community, see some
recent papers \cite{AnandCKP,Klavzar,Nowakowski,Zhu}.

In this paper, we consider four standard products: the
lexicographic, the strong, the Cartesian and the direct with respect
to the monochromatic connection number. Every of these
four products will be treated in one of the forthcoming sections. In Section $3$, we demonstrate the usefulness of the proposed
constructions by applying them to some instances of product networks.

\section{Main results}

In this section, we study the monochromatic connection number of four graph product.

\begin{lem}{\upshape \cite{CaroY}}\label{lem1}
Let $G$ be a connected graph with $n(G)$ vertices and $m(G)$ edges. Then
$$
mc(G)\leq E(G)-V(G)+\kappa(G)+1.
$$
\end{lem}

In \cite{Spacapan}, \u{S}pacapan obtained the following result.

\begin{lem}{\upshape\cite{Spacapan}}\label{lem2}
Let $G$ and $H$ be two nontrivial graphs. Then
$$
\kappa(G\Box
H)=\min\{\kappa(G)|V(H)|,\kappa(H)|V(G)|,\delta(G)+\delta(H)\}.
$$
\end{lem}

Yang and Xu \cite{Yang} investigated the classical connectivity of the lexicographic product
of two graphs.

\begin{lem}{\upshape\cite{Yang}}\label{lem3}
Let $G$ and $H$ be two graphs. If $G$ is non-trivial, non-complete
and connected, then
$$
\kappa(G\circ H)=\kappa(G)|V(H)|.
$$
\end{lem}

Let $S_G$ and $S_H$ be separating sets of connected graphs $G$ and $H$, and let $G^\prime$ and $H^\prime$ be arbitrary
connected components of $G-S_G$ and $H-S_H$. Then the set of vertices
$$
(S_G\times V(H^\prime))\cup(S_G\times S_H)\cup(V(G^\prime) \times S_H)
$$
is called \emph{a $\daleth$-set} of $G\boxtimes H$; see {\upshape\cite{Hammack}}.

\begin{lem}{\upshape\cite{Hammack}}\label{lem4}
Let $G$ and $H$ be connected graphs, at least one not complete. Set $\ell(G\boxtimes H)$ be the minimum size of a $\daleth$-set of $G\boxtimes H$. Then
$$
\kappa(G\boxtimes H)=\min\{\kappa(G)|V(H)|,\kappa(H)|V(G)|,\ell(G\boxtimes H)\}.
$$
\end{lem}

\begin{lem}{\upshape\cite{Hammack}}\label{lem5}
Let $G$ and $H$ be nonbipartite graphs. Then
$$
\kappa^\prime(G\times H)=\min\{2\kappa^\prime(G)|V(H)|,2\kappa^\prime(H)|V(G)|,\delta(G)\delta(H)\}.
$$
\end{lem}

Let $d_{G}(u,v)$ denote the distance between $u$ and $v$ in $G$.
Denote by $d_{G}(u)$ the degree of vertex $u$ in $G$. The following lemma is from {\upshape \cite{Hammack}}.

\begin{lem}{\upshape \cite{Hammack}}\label{lem6}
Let $(g,h)$ and $(g',h')$ be two vertices of $G\Box H$. Then
$$
d_{G\Box H}((g,h),(g',h'))= d_{G}(g,g')+ d_{H}(h,h').
$$
\end{lem}

\begin{cor}\label{cor1}
Let $G$ be a connected graph. Then
$$
diam(G\Box H)=diam(G)+diam(H).
$$
\end{cor}

\begin{lem}{\upshape \cite{Hammack}}\label{lem7}
Let $(g,h)$ and $(g',h')$ be two vertices of $G\circ H$. Then
$$
d_{G\circ H}((g,h),(g',h'))=\left\{
\begin{array}{ll}
d_{G}(g,g'), &if~g\neq g';\\[0.2cm]
d_{H}(h,h'), &if~g=g'~and~d_{G}(g)=0;\\[0.2cm]
\min \{d_{H}(h,h'),2\}, &if~g=g'~and~d_{G}(g)\neq 0.
\end{array}
\right.
$$
\end{lem}

\begin{lem}{\upshape \cite{Hammack}}\label{lem8}
Let $(g,h)$ and $(g',h')$ be two vertices of $G\Box H$. Then
$$
d_{G\boxtimes H}((g,h),(g',h'))=\max\{d_{G}(g,g'), d_{H}(h,h')\}.
$$
\end{lem}

\begin{cor}\label{cor2}
Let $G$ be a connected graph. Then
$$
diam(G\boxtimes H)=\max\{diam(G),diam(H)\}.
$$
\end{cor}

\subsection{The Cartesian product}

The \emph{Cartesian product} of two graphs $G$ and $H$, written as
$G\Box H$, is the graph with vertex set $V(G)\times V(H)$, in which
two vertices $(g,h)$ and $(g',h')$ are adjacent if and only if
$g=g'$ and $(h,h')\in E(H)$, or $h=h'$ and $(g,g')\in E(G)$.
Clearly, the Cartesian product is commutative, that is, $G\Box H$ is
isomorphic to $H\Box G$. The Cartesian product is commutative, that is,
$G\Box H\cong H\Box G$. Clearly, $|E(G\circ
H)|=|E(H)||V(G)|+|E(G)||V(H)|$.

\begin{thm}\label{th2}
Let $G$ and $H$ be a connected graph.

$(1)$ If neither $G$ nor $H$ is a tree, then
$$
\max\{|E(G)||V(H)|,|E(H)||V(G)|\}+2\leq mc(G\Box H)\leq |E(G)||V(H)|+(|E(H)|-1)|V(G)|+1.
$$

$(2)$ If $G$ is not a tree and $H$ is a tree, then
$$
|E(H)||V(G)|+2\leq mc(G\Box H)\leq|E(G)||V(H)|+1.
$$

$(3)$ If both $G$ and $H$ are trees, then
$$
|E(G)||E(H)|+1\leq mc(G\Box H)\leq |E(G)||E(H)|+2.
$$

Moreover, the lower bounds are sharp.
\end{thm}

\begin{pf}
$(1)$ Since $H$ is not a tree, it follows that $|E(H)|\geq |V(H)|$. By Observation \ref{obs1}, we have
\begin{eqnarray*}
mc(G\Box H)&\geq& |E(G\Box H)|-|V(G\Box H)|+2\\
&=&|E(G)||V(H)|+|E(H)||V(G)|-|V(G)||V(H)|+2\\
&\geq&|E(G)||V(H)|+2.
\end{eqnarray*}

From Lemma \ref{lem2}, $\kappa(G\Box
H)=\min\{\kappa(G)|V(H)|,\kappa(H)|V(G)|,\delta(G)+\delta(H)\}\leq \kappa(H)|V(G)|\leq (|V(H)|-1)|V(G)|$. Furthermore, by Lemma \ref{lem1}, we have
\begin{eqnarray*}
&&mc(G\Box H)\\
&\leq& E(G\Box H)-V(G\Box H)+\kappa(G\Box H)+1\\
&=&|E(G)||V(H)|+|E(H)||V(G)|-|V(G)||V(H)|+\kappa(G\Box H)+1\\
&\leq&|E(G)||V(H)|+|E(H)||V(G)|-|V(G)||V(H)|+|V(G)|(|V(H)|-1)+1\\
&=&|E(G)||V(H)|+|E(H)||V(G)|-|V(G)|+1.
\end{eqnarray*}

$(2)$ Since $G$ is not a tree and $H$ is a tree, it follows that $|E(G)|\geq |V(G)|$ and $|E(H)|=|V(H)|-1$. By Observation \ref{obs1}, we have
\begin{eqnarray*}
mc(G\Box H)&\geq& E(G\Box H)-V(G\Box H)+2\\
&=&|E(G)||V(H)|+|E(H)||V(G)|-|V(G)||V(H)|+2\\
&\geq&|V(G)||V(H)|+|E(H)||V(G)|-|V(G)||V(H)|+2\\
&=&|E(H)||V(G)|+2.
\end{eqnarray*}

From Lemma \ref{lem2}, $\kappa(G\Box
H)=\min\{\kappa(G)|V(H)|,\kappa(H)|V(G)|,\delta(G)+\delta(H)\}\leq \kappa(H)|V(G)|\leq |V(G)|$. By Lemma \ref{lem1}, we have
\begin{eqnarray*}
mc(G\Box H)&\leq& E(G\Box H)-V(G\Box H)+\kappa(G\Box H)+1\\
&=&|E(G)||V(H)|+|E(H)||V(G)|-|V(G)||V(H)|+\kappa(G\Box H)+1\\
&=&|E(G)||V(H)|+(|V(H)|-1)|V(G)|-|V(G)||V(H)|+\kappa(G\Box H)+1\\
&\leq&|E(G)||V(H)|-|V(G)|+|V(G)|+1\\
&\leq&|E(G)||V(H)|+1.
\end{eqnarray*}

$(3)$ Since both $G$ and $H$ are trees, it follows that $|E(G)|=|V(G)|-1$ and $|E(H)|=|V(H)|-1$. By Observation \ref{obs1}, we have
\begin{eqnarray*}
mc(G\Box H)&\geq& E(G\Box H)-V(G\Box H)+2\\
&=&|E(G)||V(H)|+|E(H)||V(G)|-|V(G)||V(H)|+2\\
&=&(|V(G)|-1)|V(H)|+(|V(H)|-1)|V(G)|-|V(G)||V(H)|+2\\
&=&|E(G)||E(H)|+1.
\end{eqnarray*}

From Lemma \ref{lem2}, $\kappa(G\Box
H)=\min\{\kappa(G)|V(H)|,\kappa(H)|V(G)|,\delta(G)+\delta(H)\}\leq \delta(G)+\delta(H)=2$. By Lemma \ref{lem1}, we have
\begin{eqnarray*}
mc(G\Box H)&\leq& E(G\Box H)-V(G\Box H)+\kappa(G\Box H)+1\\
&=&|E(G)||V(H)|+|E(H)||V(G)|-|V(G)||V(H)|+\kappa(G\Box H)+1\\
&=&(|V(G)|-1)|V(H)|+(|V(H)|-1)|V(G)|-|V(G)||V(H)|+\kappa(G\Box H)+1\\
&=&|E(G)||E(H)|+\kappa(G\Box H)\\
&\leq&|E(G)||E(H)|+2.
\end{eqnarray*}
\end{pf}

To show the sharpness of the lower bounds in Theorem \ref{th2}, we consider the following example.

\noindent \textbf{Example 1:} $(1)$ Let $G$ be a cycle of order at least $3$, and $H$ be a cycle of order at least $4$. From
Corollary \ref{cor1}, $diam(G\Box
H)=diam(G)+diam(H)\geq 3$. By
Theorem \ref{th1}, $mc(G\Box H)=|E(G\Box H)|-|V(G\Box H)|+2=|E(G)||V(H)|+2=|E(H)||V(G)|+2$.

$(2)$ Let $G$ be a cycle of order at least $4$, and $H$ be a path of order at least $3$. From Corollary \ref{cor1}, $diam(G\Box
H)=diam(G)+diam(H)\geq 3$. By
Theorem \ref{th1}, $mc(G\Box H)=|E(G\Box H)|-|V(G\Box H)|+2=|E(H)||V(G)|+2$.

$(3)$ Let $G=P_2$ and $H$ be a path of order at least $3$. From Corollary \ref{cor1}, $diam(G\Box
H)=diam(G)+diam(H)\geq 3$. Therefore, $mc(G\Box H)=|E(G\Box H)|-|V(G\Box H)|+2=|E(G)||E(H)|+1$.

The following corollary is immediate from Theorem \ref{th2}.

\begin{cor}\label{cor3}
Let $G$ and $H$ be a connected graph.

$(1)$ If neither $G$ nor $H$ is a tree, then
$mc(G\Box H)\geq \max\{mc(G)|V(H)|+2,mc(H)|V(G)|+2\}$.

$(2)$ If $G$ is not a tree and $H$ is a tree, then
$mc(G\Box H)\geq mc(H)|V(G)|+2$.

$(3)$ If both $G$ and $H$ are trees, then
$mc(G\Box H)\geq mc(G)mc(H)+1$.
\end{cor}

\subsection{The lexicographical product}

The lexicographic product $G\circ H$ of graphs $G$ and $H$ has the
vertex set $V(G\circ H)=V(G)\times V(H)$. Two vertices
$(g,h),(g',h')$ are adjacent if $gg'\in E(G)$, or if $g=g'$ and
$hh'\in E(H)$. The lexicographic product is not commutative and is
connected whenever $G$ is connected. Note that unlike the Cartesian Product,
the lexicographic product is a non-commutative product since $G\circ
H$ need not be isomorphic to $H\circ G$. Clearly, $|E(G\circ
H)|=|E(H)||V(G)|+|E(G)||V(H)|^2$.

\begin{thm}\label{th3}
Let $G$ and $H$ be a connected graph.

$(1)$ If neither $G$ nor $H$ is a tree, then
$$
|E(G)||V(H)|^2+2\leq mc(G\circ H)\leq |E(H)||V(G)|+|E(G)||V(H)|^2-|V(H)|+1.
$$

$(2)$ If $G$ not a tree and $H$ is a tree, then
$$
|E(H)||V(G)|(|V(H)|+1)+2\leq mc(G\circ H)\leq|E(H)||V(G)|+|E(G)||V(H)|^2-|V(H)|+1.
$$

$(3)$ If $H$ not a tree and $G$ is a tree, then
$$
|E(H)||V(G)|^2+2\leq mc(G\circ H)\leq|E(H)||V(G)|+|E(G)||V(H)|^2-|V(H)|+1.
$$

$(4)$ If both $G$ and $H$ are trees, then
$$
|E(H)||E(G)|(|V(H)|+1)+1\leq mc(G\circ H))\leq|E(H)||E(G)|(|V(H)|+1)+|V(H)|.
$$

Moreover, the lower bounds are sharp.
\end{thm}
\begin{pf}
$(1)$ Since $H$ is not a tree, it follows that $|E(H)|\geq |V(H)|$. By Observation \ref{obs1}, we have
\begin{eqnarray*}
mc(G\circ H)&\geq& |E(G\circ H)|-|V(G\circ H)|+2\\
&=&|E(H)||V(G)|+|E(G)||V(H)|^2-|V(G)||V(H)|+2\\
&\geq&|E(G)||V(H)|^2+2.
\end{eqnarray*}

From Lemma \ref{lem3}, $\kappa(G\circ
H)=\kappa(G)|V(H)|\leq (|V(G)|-1)|V(H)|$. By Lemma \ref{lem1}, we have
\begin{eqnarray*}
mc(G\circ H)&\leq& |E(G\circ H)|-|V(G\circ H)|+\kappa(G\circ H)+1\\
&=&|E(H)||V(G)|+|E(G)||V(H)|^2-|V(G)||V(H)|+\kappa(G\circ H)+1\\
&\leq&|E(H)||V(G)|+|E(G)||V(H)|^2-|V(G)||V(H)|+(|V(G)|-1)|V(H)|+1\\
&=&|E(H)||V(G)|+|E(G)||V(H)|^2-|V(H)|+1.
\end{eqnarray*}

$(2)$ Since $G$ is not a tree and $H$ is a tree, it follows that $|E(G)|\geq |V(G)|$ and $|E(H)|=|V(H)|-1$. By Observation \ref{obs1}, we have
\begin{eqnarray*}
mc(G\circ H)&\geq& |E(G\circ H)|-|V(G\circ H)|+2\\
&=&|E(H)||V(G)|+|E(G)||V(H)|^2-|V(G)||V(H)|+2\\
&=&(|V(H)|-1)|V(G)|+|E(G)||V(H)|^2-|V(G)||V(H)|+2\\
&\geq&|E(H)||V(G)|(|V(H)|+1)+2.
\end{eqnarray*}

From Lemma \ref{lem3}, $\kappa(G\circ
H)=\kappa(G)|V(H)|\leq (|V(G)|-1)|V(H)|$. By Lemma \ref{lem1}, we have
\begin{eqnarray*}
&&mc(G\circ H)\\
&\leq& |E(G\circ H)|-|V(G\circ H)|+\kappa(G\circ H)+1\\
&=&|E(H)||V(G)|+|E(G)||V(H)|^2-|V(G)||V(H)|+\kappa(G\circ H)+1\\
&=&|E(H)||V(G)|+|E(G)||V(H)|^2-|V(G)||V(H)|+\kappa(G\circ H)+1\\
&\leq&|E(H)||V(G)|+|E(G)||V(H)|^2-|V(G)||V(H)|+(|V(G)|-1)|V(H)|+1\\
&=&|E(H)||V(G)|+|E(G)||V(H)|^2-|V(H)|+1.
\end{eqnarray*}

$(3)$ Since $H$ is not a tree, it follows that $|E(H)|\geq |V(H)|$. By Observation \ref{obs1}, we have
\begin{eqnarray*}
mc(G\circ H)&\geq& |E(G\circ H)|-|V(G\circ H)|+2\\
&=&|E(H)||V(G)|+|E(G)||V(H)|^2-|V(G)||V(H)|+2\\
&\geq&|E(G)||V(H)|^2+2.
\end{eqnarray*}

From Lemma \ref{lem3}, $\kappa(G\circ
H)=\kappa(G)|V(H)|=|V(H)|$. By Lemma \ref{lem1}, we have
\begin{eqnarray*}
mc(G\circ H)&\leq& |E(G\circ H)|-|V(G\circ H)|+\kappa(G\circ H)+1\\
&=&|E(H)||V(G)|+|E(G)||V(H)|^2-|V(G)||V(H)|+\kappa(G\circ H)+1\\
&=&|E(H)||V(G)|+|E(G)||V(H)|^2-|V(G)||V(H)|+\kappa(G\circ H)+1\\
&=&|E(H)||V(G)|+|E(G)||V(H)|^2-|V(G)||V(H)|+|V(H)|+1.
\end{eqnarray*}

$(4)$ Since both $G$ and $H$ are trees, it follows that $|E(G)|=|V(G)|-1$ and $|E(H)|=|V(H)|-1$. By Observation \ref{obs1}, we have
\begin{eqnarray*}
mc(G\circ H)&\geq& |E(G\circ H)|-|V(G\circ H)|+2\\
&=&|E(H)||V(G)|+|E(G)||V(H)|^2-|V(G)||V(H)|+2\\
&=&(|V(H)|-1)|V(G)|+(|V(G)|-1)|V(H)|^2-|V(G)||V(H)|+2\\
&=&|E(H)||E(G)|(|V(H)|+1)+1.
\end{eqnarray*}

From Lemma \ref{lem3}, $\kappa(G\circ
H)=\kappa(G)|V(H)|=|V(H)|$. By Lemma \ref{lem1}, we have
\begin{eqnarray*}
&&mc(G\circ H)\\
&\leq& |E(G\circ H)|-|V(G\circ H)|+\kappa(G\circ H)+1\\
&=&|E(H)||V(G)|+|E(G)||V(H)|^2-|V(G)||V(H)|+\kappa(G\circ H)+1\\
&=&(|V(H)|-1)|V(G)|+(|V(G)|-1)|V(H)|^2-|V(G)||V(H)|+\kappa(G\circ H)+1\\
&=&|E(H)||E(G)|(|V(H)|+1)+\kappa(G\circ H)\\
&=&|E(H)||E(G)|(|V(H)|+1)+V(H).
\end{eqnarray*}
\end{pf}

To show the sharpness of the lower bounds in Theorem \ref{thm3}, we consider the following example.

\noindent \textbf{Example 2:} $(1)$ Let $G$ be a cycle of order at least $6$, and $H$ be a cycle of order at least $3$.
From Lemma \ref{lem7}, $diam(G\circ H)\geq diam(G)\geq 3$. Therefore, $mc(G\circ H)=|E(G\circ H)|-|V(G\circ H)|+2=|E(G)||V(H)|^2+2$.

$(2)$ Let $G$ be a cycle of order at least $6$, and $H=P_n, n\geq 4$. By
Lemma \ref{lem7}, $diam(G\circ H)\geq diam(G)\geq 3$. Therefore, $mc(G\circ H)=|E(H)||V(G)|(|V(H)|+1)+2$.

$(3)$ Let $G$ be a path of order at least $4$, and $H$ be a cycle of order at least $3$. By
Lemma \ref{lem7}, $diam(G\circ H)\geq diam(G)\geq 3$. Therefore, $mc(G\circ H)=|E(H)||V(G)|^2+2$.

$(4)$ Let $G$ be a path of order at least $4$, and $H=P_2$. By
Lemma \ref{lem7}, $diam(G\circ H)\geq diam(G)\geq 3$. Therefore, $mc(G\circ H)=|E(H)||E(G)|(|V(H)|+1)+1$.

The following corollary is immediate from Theorem \ref{th3}.

\begin{cor}\label{cor4}
Let $G$ and $H$ be a connected graph.

$(1)$ If neither $G$ nor $H$ is a tree, then
$mc(G\circ H)\geq mc(G)|V(H)|^2+2$.

$(2)$ If $G$ not a tree and $H$ is a tree, then
$mc(G\circ H)\geq mc(H)|V(G)|(|V(H)|+1)+2$.

$(3)$ If $H$ not a tree and $G$ is a tree, then
$mc(G\circ H)\geq mc(H)|V(G)|^2+2$.

$(4)$ If both $G$ and $H$ are trees, then
$mc(G\circ H))\geq mc(G)mc(H)(|V(H)|+1)+1$.

Moreover, the lower bounds are sharp.
\end{cor}

\subsection{The strong product}

The \emph{strong product} $G\boxtimes H$ of graphs $G$ and $H$ has
the vertex set $V(G)\times V(H)$. Two vertices $(g,h)$ and $(g',h')$
are adjacent whenever $gg'\in E(G)$ and $h=h'$, or $g=g'$ and $hh'
\in E(H)$, or $gg'\in E(G)$ and $hh'\in E(H)$. Clearly, $|E(G\boxtimes
H)|=|E(H)||V(G)|+|E(G)||V(H)|+2|E(G)||E(H)|$.

\begin{thm}\label{thm4}
Let $G$ and $H$ be a connected graph, and at least one of $G$ and $H$ is not a complete graph.

$(1)$ If neither $G$ nor $H$ is a tree, then
$$
mc(G\boxtimes H)\geq \max\{|E(G)||V(H)|+2|E(H)||E(G)|+2,|E(H)||V(G)|+2|E(H)||E(G)|+2\}
$$
and
$$
mc(G\boxtimes H)\leq |E(G)||V(H)|+|E(H)||V(G)|+2|E(H)||E(G)|-\min\{|V(H)|,|V(G)|\}+1.
$$

$(2)$ If $G$ not a tree and $H$ is a tree, then
$$
|E(H)||V(G)|+2|E(H)||E(G)|+2\leq mc(G\boxtimes H)\leq |E(G)||V(H)|+2|E(H)||E(G)|+1.
$$

$(3)$ If both $G$ and $H$ are trees, then
$$
3|E(H)||E(G)|+1\leq mc(G\boxtimes H)\leq 3|E(H)||E(G)|+\min\{|V(G)|,|V(H)|\}.
$$
Moreover, the lower bounds are sharp.
\end{thm}

\begin{pf}
$(1)$ Since $H$ is not a tree, it follows that $|E(H)|\geq |V(H)|$. By Observation \ref{obs1}, we have
\begin{eqnarray*}
&&mc(G\boxtimes H)\\
&\geq& |E(G\boxtimes H)|-|V(G\boxtimes H)|+2\\
&=&|E(G)||V(H)|+|E(H)||V(G)|+2|E(H)||E(G)|-|V(G)||V(H)|+2\\
&\geq&|E(G)||V(H)|+2|E(H)||E(G)|+2.
\end{eqnarray*}

From Lemma \ref{lem4}, $\kappa(G\boxtimes
H)=\min\{\kappa(G)|V(H)|,\kappa(H)|V(G)|,\ell(G\boxtimes H)\}\leq \min\{(|V(G)|-1)|V(H)|,(|V(H)|-1)|V(G)|\}=|V(G)||V(H)|-\min\{|V(H)|,|V(G)|\}$. By Lemma \ref{lem1}, we have
\begin{eqnarray*}
&&mc(G\boxtimes H)\\
&\leq& E(G\boxtimes H)-V(G\boxtimes H)+\kappa(G\boxtimes H)+1\\
&=&|E(G)||V(H)|+|E(H)||V(G)|+2|E(H)||E(G)|-|V(G)||V(H)|+\kappa(G\boxtimes H)+1\\
&\leq&|E(G)||V(H)|+|E(H)||V(G)|+2|E(H)||E(G)|-\min\{|V(H)|,|V(G)|\}+1.
\end{eqnarray*}

$(2)$ Since $G$ is not a tree, it follows that $|E(G)|\geq |V(G)|$. Since $H$ is a tree, we have $|E(H)|=|V(H)|-1$. By Observation \ref{obs1}, we have
\begin{eqnarray*}
&&mc(G\boxtimes H)\\
&\geq& |E(G\boxtimes H)|-|V(G\boxtimes H)|+2\\
&=&|E(G)||V(H)|+|E(H)||V(G)|+2|E(H)||E(G)|-|V(G)||V(H)|+2\\
&=&|E(G)||V(H)|+(|V(H)|-1)|V(G)|+2|E(H)||E(G)|-|V(G)||V(H)|+2\\
&\geq&|E(H)||V(G)|+2|E(H)||E(G)|+2.
\end{eqnarray*}

From Lemma \ref{lem4}, $\kappa(G\boxtimes
H)=\min\{\kappa(G)|V(H)|,\kappa(H)|V(G)|,\ell(G\boxtimes H)\}\leq \kappa(H)|V(G)|=|V(G)|$. By Lemma \ref{lem1}, we have
\begin{eqnarray*}
&&mc(G\boxtimes H)\\
&\leq& |E(G\boxtimes H)|-|V(G\boxtimes H)|+\kappa(G\boxtimes H)+1\\
&=&|E(G)||V(H)|+|E(H)||V(G)|+2|E(H)||E(G)|-|V(G)||V(H)|+\kappa(G\boxtimes H)+1\\
&=&|E(G)||V(H)|+(|V(H)|-1)|V(G)|+2|E(H)||E(G)|-|V(G)||V(H)|+\kappa(G\boxtimes H)+1\\
&=&|E(G)||V(H)|-|V(G)|+2|E(H)||E(G)|+\kappa(G\boxtimes H)+1\\
&\leq&|E(G)||V(H)|+2|E(H)||E(G)|+1.
\end{eqnarray*}

$(3)$ Since both $G$ and $H$ are trees, it follows that $|E(G)|=|V(G)|-1$ and $|E(H)|=|V(H)|-1$. By Observation \ref{obs1}, we have
\begin{eqnarray*}
&&mc(G\boxtimes H)\\
&\geq& |E(G\boxtimes H)|-|V(G\boxtimes H)|+2\\
&=&|E(G)||V(H)|+|E(H)||V(G)|+2|E(H)||E(G)|-|V(G)||V(H)|+2\\
&=&(|V(G)|-1)|V(H)|+(|V(H)|-1)|V(G)|+2|E(H)||E(G)|-|V(G)||V(H)|+2\\
&=&3|E(H)||E(G)|+1.
\end{eqnarray*}

From Lemma \ref{lem4}, $\kappa(G\boxtimes
H)=\min\{\kappa(G)|V(H)|,\kappa(H)|V(G)|,\ell(G\boxtimes H)\}\leq \min\{\kappa(G)|V(H)|,$ $\kappa(H)|V(G)|\}\leq \min\{|V(H)|,|V(G)|\}$. By Lemma \ref{lem1}, we have
\begin{eqnarray*}
&&mc(G\boxtimes H)\\
&\leq& |E(G\boxtimes H)|-|V(G\boxtimes H)|+\kappa(G\boxtimes H)+1\\
&=&|E(G)||V(H)|+|E(H)||V(G)|+2|E(H)||E(G)|-|V(G)||V(H)|+\kappa(G\boxtimes H)+1\\
&=&|V(G)|-1)|V(H)|+(|V(H)|-1)|V(G)|+2|E(H)||E(G)|-|V(G)||V(H)|+\kappa(G\boxtimes H)+1\\
&=&3|E(H)||E(G)|+\kappa(G\boxtimes H)\\
&\leq&3|E(H)||E(G)|+\min\{|V(H)|,|V(G)|\}).
\end{eqnarray*}
\end{pf}

To show the sharpness of the lower bounds in Theorem \ref{thm4}, we consider the following example.

\noindent \textbf{Example 3:} $(1)$ Let $G$ be a cycle of order at least $6$, and $H$ be a cycle of order at least $3$. By
Corollary \ref{cor2}, $diam(G\boxtimes H)=\max\{diam(G),diam(H)\}\geq 3$. Therefore, $mc(G\boxtimes H)=|E(G)||V(H)|+2|E(H)||E(G)|+2=|E(H)||V(G)|+2|E(H)||E(G)|+2$.

$(2)$ Let $G$ be a cycle of order at least $3$, and $H$ be a cycle of order at least $4$. By
Corollary \ref{cor2}, $diam(G\boxtimes H)=\max\{diam(G),diam(H)\}\geq 3$. Therefore, $mc(G\boxtimes H)=|E(H)||V(G)|+2|E(H)||E(G)|+2$.

$(3)$ Let $G=P_2$ and $H$ be a cycle of order at least $4$. By
Corollary \ref{cor2}, $diam(G\boxtimes H)=\max\{diam(G),diam(H)\}\geq 3$. Therefore, $mc(G\boxtimes H)=3|E(H)||E(G)|+1$.

The following corollary is immediate from Theorem \ref{thm4}.

\begin{cor}\label{cor5}
Let $G$ and $H$ be a connected graph.

$(1)$ If neither $G$ nor $H$ is a tree, then
$$
mc(G\boxtimes H)\geq \max\{|mc(G)||V(H)|+2|mc(H)||mc(G)|+2,|mc(H)||V(G)|+2|mc(H)||mc(G)|+2\}.
$$

$(2)$ If $G$ not a tree and $H$ is a tree, then
$$
mc(G\boxtimes H)\geq |mc(H)||V(G)|+2|mc(H)||mc(G)|+2.
$$

$(3)$ If both $G$ and $H$ are trees, then
$$
mc(G\boxtimes H)\geq 3|mc(H)||mc(G)|+1.
$$
Moreover, the lower bounds are sharp.
\end{cor}

\subsection{The direct product}

The \emph{direct product} $G\times H$ of graphs $G$ and $H$ has the vertex
set $V(G)\times V(H)$. Two vertices $(g,h)$ and $(g',h')$ are
adjacent if the projections on both coordinates are adjacent, i.e.,
$gg'\in E(G)$ and $hh'\in E(H)$. Clearly, $|E(G\times
H)|=2|E(G)||E(H)|$.

\begin{thm}\label{thm5}
Let $G$ and $H$ be nonbipartite graphs. Then
$$
|E(H)||E(G)|+2\leq mc(G\times H)\leq 2|E(H)||E(G)|+1.
$$

Moreover, the lower bounds are sharp.
\end{thm}
\begin{pf}
Since $H$ is not a tree, it follows that $|E(H)|\geq |V(H)|$. By Observation \ref{obs1}, we have
\begin{eqnarray*}
mc(G\times H)&\geq& |E(G\times H)|-|V(G\times H)|+2\\
&=&2|E(H)||E(G)|-|V(G)||V(H)|+2\\
&\geq&|E(H)||E(G)|+2.
\end{eqnarray*}

From Lemma \ref{lem5}, $\kappa(G\times
H)\leq \kappa^\prime(G\times H)=\min\{2\kappa^\prime(G)|V(H)|,2\kappa^\prime(H)|V(G)|,\delta(G)\delta(H)\}\leq \delta(G)\delta(H)\leq |V(G)||V(H)|$. By Lemma \ref{lem1}, we have
\begin{eqnarray*}
mc(G\times H)&\leq& E(G\times H)-V(G\times H)+\kappa(G\times H)+1\\
&=&2|E(H)||E(G)|-|V(G)||V(H)|+\kappa(G\times H)+1\\
&=&2|E(H)||E(G)|+1.
\end{eqnarray*}
\end{pf}

To show the sharpness of the lower bounds in Theorem \ref{thm5}, we consider the following example.

\noindent \textbf{Example 4:} Let $G$ be a cycle of order at least $3$, and $H$ be a cycle of order at least $6$. By
Lemma \ref{lem2}, $diam(G\times H)=\max\{diam(G),diam(H)\}\geq 3$. Therefore, $mc(G\times H)=|E(H)||E(G)|+2$.

The following corollary is immediate from Theorem \ref{thm5}.

\begin{cor}\label{cor4}
Let one of $G$ and $H$ be a non-bipartite connected graph. Then
$$
mc(G\times H)\geq |mc(H)||mc(G)|+2.
$$
\end{cor}

\section{Applications}

In this section, we demonstrate the usefulness of the proposed
constructions by applying them to some instances of Cartesian and
lexicographical product networks.

\subsection{Two-dimensional grid graph}

A \emph{two-dimensional grid graph} is an $m\times n$ graph
$G_{n,m}$ that is the graph Cartesian product $P_n\Box P_m$ of path
graphs on $m$ and $n$ vertices. See Figure 1 $(a)$ for the case
$m=3$. For more details on grid graph, we refer to \cite{Calkin,
Itai}. The network $P_n\circ P_m$ is the graph lexicographical
product $P_n\circ P_m$ of path graphs on $m$ and $n$ vertices. For
more details on $P_n\circ P_m$, we refer to \cite{Mao2}. See Figure
1 $(b)$ for the case $m=3$.

\begin{figure}[!hbpt]
\begin{center}
\includegraphics[scale=0.9]{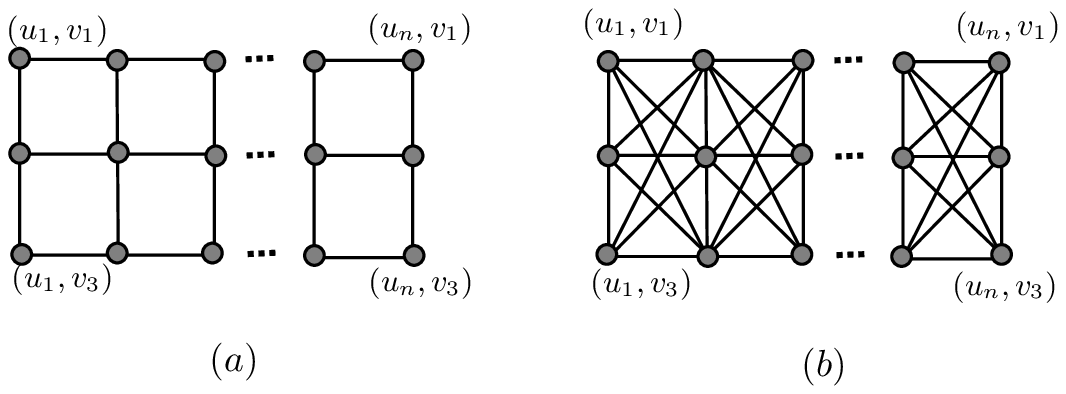}\\
Figure 1: $(a)$ Two-dimensional grid graph $G_{n,3}$; $(b)$ The
network $P_n\circ P_3$.
\end{center}\label{fig7}
\end{figure}

\begin{pro}\label{pro1}
$(i)$ For network $P_n\Box P_m \ (n\geq 3, m\geq 2)$,
$$
mc(P_n\Box P_m)=nm-n-m+2.
$$

$(ii)$ For network $P_n\circ P_m \ (n\geq 4, m\geq 3)$,
$$
mc(P_n\circ P_m)=m^2n-m^2-n+2.
$$
\end{pro}
\begin{pf}
$(i)$ From Corollary \ref{cor1}, $diam(G\Box
H)=diam(G)+diam(H)\geq 3$. Therefore, From Theorem \ref{th1}, we have
\begin{eqnarray*}
mc(P_n\Box P_m)
&=&|E(P_n\Box P_m)|-|V(P_n\Box P_m)|+2\\
&=&(|V(P_n)|-1)|V(P_m)|+(|V(P_m)|-1)|V(P_n)|-|V(P_n)||V(P_m)|+2\\
&=&(n-1)m+(m-1)n-nm+2\\
&=&nm-n-m+2.
\end{eqnarray*}

$(2)$ From Lemma \ref{lem7}, $diam(G\circ H)\geq diam(G)\geq 3$. From Theorem \ref{th1}, we have
\begin{eqnarray*}
mc(P_n\circ P_m)&=& E(P_n\circ P_m)-V(P_n\circ P_m)+2\\
&=&|E(P_m)||V(P_n)|+|E(P_n)||V(P_m)|^2-|V(P_n)||V(P_m)|+2\\
&=&(m-1)n+(n-1)m^2-mn+2\\
&=&m^2n-m^2-n+2.
\end{eqnarray*}
\end{pf}

\subsection{$n$-dimensional mesh}

An \emph{$n$-dimensional mesh} is the Cartesian product of $n$
linear arrays. By this definition, two-dimensional grid graph is a
$2$-dimensional mesh. An \emph{$n$-dimensional hypercube} is a special case
of an $n$-dimensional mesh, in which the $n$ linear arrays are all
of size $2$; see \cite{Johnsson}.

\begin{pro}\label{pro2}
$(i)$ For $n$-dimensional mesh $P_{L_1}\Box P_{L_2}\Box \cdots \Box
P_{L_n} \ (n\geq 4)$,
$$
mc(P_{L_1}\Box P_{L_2}\Box \cdots \Box P_{L_n})\geq (2\ell_1\ell_2-\ell_1-\ell_2)(\ell_3\ell_4\cdots \ell_n)+2.
$$

$(ii)$ For network $P_{L_1}\circ P_{L_2}\circ
\cdots \circ P_{L_n}$,
$$mc(P_{L_1}\circ P_{L_2}\circ \cdots \circ P_{L_n})\geq (\ell_1\ell^2_2+\ell_1\ell_2-\ell_1-\ell^2_2)(\ell_3\ell_4\cdots \ell_n)^2+2.$$
\end{pro}
\begin{pf}
$(i)$ By Lemma \ref{lem2}, we have $diam((P_{L_1}\Box P_{L_2}\Box
\cdots \Box P_{L_n})=\sum_{i=1}^ndiam(P_{L_i})\geq 3$. Set $G=P_{L_1}\Box P_{L_2}$ and $H=P_{L_2}\Box \cdots \Box P_{L_n}$. Since both $G$ and $H$ are not trees, it follows from Theorem \ref{th2} that $mc(G\Box H)\geq \max\{|E(G)||V(H)|,|E(H)||V(G)|\}+2\geq |E(G)||V(H)|+2$.
From Theorem \ref{th1}, we have
\begin{eqnarray*}
&&mc(P_{L_1}\Box P_{L_2}\Box \cdots \Box P_{L_n})\\
&\geq& |E(P_{L_1}\Box P_{L_2})||V(P_{L_3}\Box \cdots \Box P_{L_n})|+2\\
&=&|E(P_{L_1})||V(P_{L_2})|+|E(P_{L_2})||V(P_{L_1})|)(|V(P_{L_3})|\cdots|V(P_{L_n})|)+2\\
&=&(2\ell_1\ell_2-\ell_1-\ell_2)\ell_3\ell_4\cdots \ell_n+2.
\end{eqnarray*}

$(ii)$ By Lemma \ref{lem7}, we have $diam((P_{L_1}\circ P_{L_2}\circ
\cdots \circ P_{L_n})=\max \{diam(R_i)\}\geq 3$. Set $G=P_{L_1}\circ P_{L_2}$ and $H=P_{L_2}\circ \cdots \circ P_{L_n}$. Since both $G$ and $H$ are not trees, it follows from Theorem \ref{th3} that $mc(G\circ H)\geq |E(G)||V(H)|^2+2$.
From Theorem \ref{th1}, we have
\begin{eqnarray*}
&&mc(P_{L_1}\circ P_{L_2}\circ \cdots \circ P_{L_n})\\
&\geq& |E(P_{L_1}\circ P_{L_2})||V(P_{L_3}\circ \cdots \circ P_{L_n})|^2+2\\
&=&(|E(P_{L_2})||V(P_{L_1})|+|E(P_{L_1})||V(P_{L_2})|^2)(|V(P_{L_3})|\cdots|V(P_{L_n})|)^2+2\\
&=&(\ell_1\ell^2_2+\ell_1\ell_2-\ell_1-\ell^2_2)(\ell_3\ell_4\cdots \ell_n)^2+2.
\end{eqnarray*}\qed
\end{pf}

\subsection{$n$-dimensional torus}

An \emph{$n$-dimensional torus} is the Cartesian product of $n$
rings $R_1,R_2,\cdots,R_n$ of size at least three. (A
ring is a cycle in Graph Theory.) The rings $R_i$ are not necessary
to have the same size. Here, we consider the
networks constructed by $R_1\Box R_2\Box \cdots \Box R_n$ and
$R_1\circ R_2\circ \cdots \circ R_n$.

\begin{pro}\label{pro3}
$(i)$ For network $R_1\Box R_2\Box \cdots \Box R_n$, $n\geq 4$
$$
mc(R_1\Box R_2\Box \cdots \Box R_n)\geq r_1r_2\cdots r_n+2.
$$
where $r_i$ is the order of $R_i$ and $3\leq i\leq
n$.

$(ii)$ For network $R_1\circ R_2\circ \cdots \circ R_n$, $n\geq 4$
$$
mc(R_1\circ R_2\circ \cdots \circ R_n)\geq r_1(r_2\cdots r_n)^2+2.
$$
\end{pro}
\begin{pf}
$(i)$ By Lemma \ref{lem2}, we have $diam((R_1\Box R_2\Box
\cdots \Box R_n)=\sum_{i=1}^ndiam(R_i)\geq 3$. Set $G=R_1$ and $H=R_2\Box
\cdots \Box R_n$. Since both $G$ and $H$ are not trees, it follows from Theorem \ref{th2} that $mc(G\Box H)\geq \max\{|E(G)||V(H)|,|E(H)||V(G)|\}+2\geq |E(G)||V(H)|+2$.
From Theorem \ref{th1}, we have
\begin{eqnarray*}
mc(R_1\Box R_2\Box \cdots \Box R_n)
&\geq& |E(R_1)||V(R_2\Box \cdots \Box R_n)|+2\\
&=&r_1r_2\cdots r_n+2.
\end{eqnarray*}

$(ii)$ By Lemma \ref{lem7}, we have $diam((R_1\circ R_2\circ \cdots \circ R_n)=\max \{diam(R_i)\}\geq 3$. Set $G=R_1$ and $H=R_2\circ \cdots \circ R_n$. Since both $G$ and $H$ are not trees, it follows from Theorem \ref{th3} that $mc(G\circ H)\geq |E(G)||V(H)|^2+2$.
From Theorem \ref{th1}, we have
\begin{eqnarray*}
mc(R_1\circ R_2\circ \cdots \circ R_n)
&\geq& |E(R_1)||V(R_1\circ R_2\circ \cdots \circ R_n)|^2+2\\
&=&r_1(r_2\cdots r_n)^2+2.
\end{eqnarray*}\qed
\end{pf}

\subsection{$n$-dimensional generalized hypercube}

Let $K_m$ be a clique of $m$ vertices, $m\geq 2$. An
\emph{$n$-dimensional generalized hypercube} \cite{DayA,
Fragopoulou} is the Cartesian product of $m$ cliques. We have the following:

\begin{pro}\label{pro4}
$(i)$ For network $K_{m_1}\Box K_{m_2}\Box \cdots \Box K_{m_n} \
(m_i\geq 2, \ n\geq 3, \ 1\leq i\leq n)$
$$
mc(K_{m_1}\Box K_{m_2}\Box \cdots \Box K_{m_n})\geq {{m_1}\choose{2}}m_2\cdots m_n+2.
$$

$(ii)$ For network $K_{m_1}\circ K_{m_2}\circ \cdots \circ K_{m_n}$,
$$
mc(K_{m_1}\circ K_{m_2}\circ \cdots \circ K_{m_n})={{m_1m_2\cdots m_n}\choose{2}}.
$$
\end{pro}
\begin{pf}
$(i)$ By Lemma \ref{lem2}, we have $diam(K_{m_1}\Box K_{m_2}\Box \cdots \Box K_{m_n})=\sum_{i=1}^nK_{m_i}\geq 3$. Set $G=K_{m_1}$ and $H=K_{m_2}\Box \cdots \Box K_{m_n}$. Since both $G$ and $H$ are not trees, it follows from Theorem \ref{th2} that $mc(G\Box H)\geq \max\{|E(G)||V(H)|,|E(H)||V(G)|\}+2\geq |E(G)||V(H)|+2$.
From Theorem \ref{th1}, we have
\begin{eqnarray*}
mc(K_{m_1}\Box K_{m_2}\Box \cdots \Box K_{m_n})
&\geq& |E(K_{m_1})||V(K_{m_2}\Box K_{m_3}\Box \cdots \Box K_{m_n})|+2\\
&=&{{m_1}\choose{2}}m_2\cdots m_n+2.
\end{eqnarray*}

$(ii)$
Note that $K_{m_1}\circ K_{m_2}\circ \cdots \circ K_{m_n}$ is a complete graph of order $\prod_{i=1}^nm_i$. From Theorem \ref{th1}, we have
\begin{eqnarray*}
mc(K_{m_1}\circ K_{m_2}\circ \cdots \circ K_{m_n})
&=&{{m_1m_2\cdots m_n}\choose{2}}.
\end{eqnarray*}\qed
\end{pf}

\subsection{$n$-dimensional hyper Petersen network}

An \emph{$n$-dimensional hyper Petersen network}
$HP_n$ is the Cartesian product of $Q_{n-3}$ and
the well-known Petersen graph \cite{Das}, where $n\geq 3$ and
$Q_{n-3}$ denotes an $(n-3)$-dimensional hypercube. The cases $n=3$
and $4$ of hyper Petersen networks are depicted in Figure 2. Note
that $HP_3$ is just the Petersen graph (see Figure 2 $(a)$).

The network $HL_n$ is the lexicographical product of
$Q_{n-3}$ and the Petersen graph, where $n\geq 3$ and $Q_{n-3}$
denotes an $(n-3)$-dimensional hypercube; see \cite{Mao2}. Note that
$HL_3$ is just the Petersen graph, and $HL_4$ is a graph obtained
from two copies of the Petersen graph by add one edge between one
vertex in a copy of the Petersen graph and one vertex in another
copy. See Figure 2 $(c)$ for an example (We only show the edges
$v_1u_i \ (1\leq i\leq 10)$).

\begin{figure}[!hbpt]
\begin{center}
\includegraphics[scale=0.8]{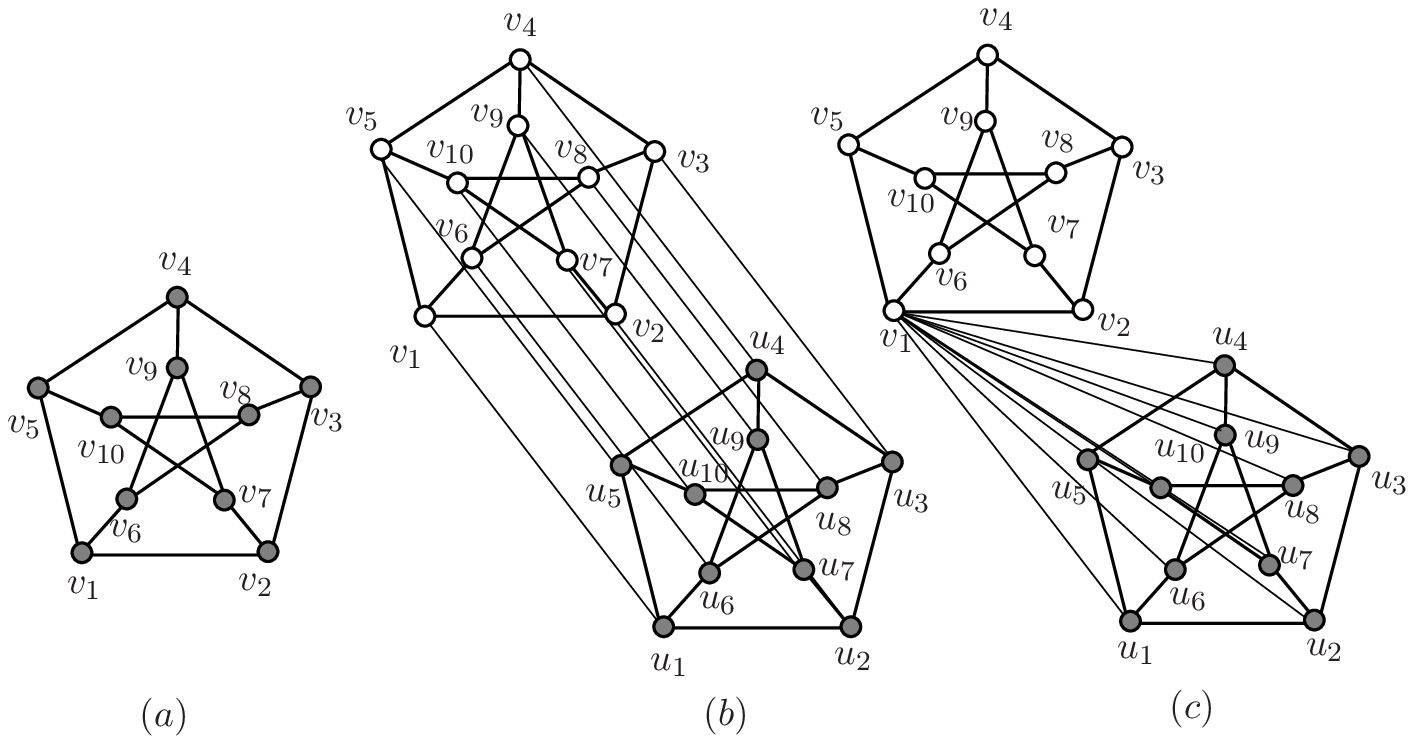}\\
Figure 2: $(a)$ Petersen graph; $(b)$ The network $HP_4$; $(c)$ The
structure of $HL_4$.
\end{center}\label{fig7}
\end{figure}

\begin{pro}\label{pro5}
$(1)$ For network $HP_3$ and $HL_3$, $mc(HP_3)=mc(HL_3)=7$;

$(2)$ For network $HL_4$ and $HP_4$, $mc(HP_4)=22$ and $112\leq mc(HL_4)\leq 121$.
\end{pro}
\begin{pf}
$(1)$ By Theorem \ref{th1}, we have $mc(HP_3)=mc(HL_3)=|E(HL_3)|-|V(HL_3)|+2=7$.

$(2)$ By Lemma \ref{lem1}, we have $mc(HP_4)=|E(HL_4)|-|V(HL_4)|+2=22$ and
$121 \geq E(HL_4)-V(HL_4)+\kappa(HL_4)+1 \geq mc(HL_4)\geq E(HL_4)-V(HL_4)+2=112$.\qed
\end{pf}

\end{document}